\documentclass[12pt]{amsart}
\usepackage{amsmath}

\usepackage{color}	 
\usepackage{soul} 

\usepackage{latexsym}
\usepackage{amsfonts, amssymb}

\newenvironment{thm}{\subsection{}{\textbf {Theorem.}}\em}{\bigskip}
\newenvironment{prop}{\subsection{}{\textbf {Proposition.}}\em}{\bigskip}
\newenvironment{cor}{\subsection{}{\textbf {Corollary.}}\em}{\bigskip}
\newenvironment{lem}{\subsection{}{\textbf {Lemma.}}\em}{\bigskip}

\newenvironment{eg}{\subsection{}{\textbf {Example.}}}{\smallskip}

\newenvironment{rem}{\subsection{}{\textbf {Remark.}}}{\smallskip}

\newcommand\fA{\ensuremath{\mathfrak A}}

\newcommand\fH{\ensuremath{\mathfrak H}}

\newcommand\fK{\ensuremath{\mathfrak K}}
\newcommand\fL{\ensuremath{\mathfrak L}}
\newcommand\fM{\ensuremath{\mathfrak M}}

\newcommand\fS{\ensuremath{\mathfrak S}}

\newcommand\fW{\ensuremath{\mathfrak W}}

\newcommand\cA{\ensuremath{\mathcal A}}
\newcommand\cB{\ensuremath{\mathcal B}}
\newcommand\cC{\ensuremath{\mathcal C}}
\newcommand\cD{\ensuremath{\mathcal D}}
\newcommand\cE{\ensuremath{\mathcal E}}

\newcommand\cH{\ensuremath{\mathcal H}}

\newcommand\cK{\ensuremath{\mathcal K}}

\newcommand\cM{\ensuremath{\mathcal M}}
\newcommand\cN{\ensuremath{\mathcal N}}
\newcommand\cO{\ensuremath{\mathcal O}}
\newcommand\cP{\ensuremath{\mathcal P}}

\newcommand\cR{\ensuremath{\mathcal R}}

\newcommand\cT{\ensuremath{\mathcal T}}
\newcommand\cU{\ensuremath{\mathcal U}}
\newcommand\cV{\ensuremath{\mathcal V}}

\newcommand\bbC{\ensuremath{\mathbb C}}
\newcommand\bbD{\ensuremath{\mathbb D}}

\newcommand\bbF{\ensuremath{\mathbb F}}

\newcommand\bbM{\ensuremath{\mathbb M}}
\newcommand\bbN{\ensuremath{\mathbb N}}

\newcommand\bbZ{\ensuremath{\mathbb Z}}

\newcommand\hilb{\ensuremath{\mathcal H}}

\newcommand\eop{{{\hfil \ensuremath \Box}}}
\newcommand\eps{\ensuremath {\varepsilon}}
\newcommand\norm{\ensuremath {\Vert}}

\newcommand\bofh{\ensuremath{\cB ( \cH)}}

\newcommand\ran{\ensuremath {\mathrm {ran}}}

\newcommand\sumin{\ensuremath {\sum_{i=1}^{n}}}

\definecolor{lightGreen}{rgb}{.75, 1, .75} 
 
\definecolor{ForestGreen}   {cmyk}{0.91,0,0.88,0.12}

\newcommand{\ltsr}{\mathrm{ltsr}}
\newcommand{\rtsr}{\mathrm{rtsr}}
\newcommand{\tsr}{\mathrm{tsr}}

\newcommand{\lbsr}{\mathrm{lBsr}}
\newcommand{\rbsr}{\mathrm{rBsr}}
\newcommand{\bsr}{\mathrm{Bsr}}

\newcommand{\qand}{\quad\text{and}\quad}
\newcommand{\rank}{\operatorname{rank}}
\newenvironment{pf}{\noindent{\textbf {Proof.}}}
{\hfill\eop\bigskip}

\begin{document}
\title{On the topological stable rank of  non-selfadjoint operator algebras}

\author
	[K.R. Davidson]{{K.R. Davidson${}^1$}}
\email{krdavids@uwaterloo.ca \vspace{-1.5ex}}	
\address
	{Pure Math.\ Dept.\\
	University of Waterloo\\
	Waterloo, Ontario\\
	Canada  \  N2L 3G1}		
\author
	[R. Levene]{{R. Levene}}	
\email{r.levene@qub.ca.uk \vspace{-1.5ex}}	
\address
	{School of Mathematics and Physics\\
	Queen's University Belfast\\
	Belfast BT7 1NN, Northern Ireland, UK }
\author
	[L.W. Marcoux]{{L.W.~Marcoux${}^1$}}
\email{LWMarcoux@uwaterloo.ca \vspace{-1.5ex}}	
\address
	{Pure Math.\ Dept.\\
	University of Waterloo\\
	Waterloo, Ontario\\
	Canada  \  N2L 3G1}	
\author
	[H. Radjavi]{{H.~Radjavi${}^1$}}	
\email{hradjavi@uwaterloo.ca \vspace{-1.5ex}}
\address
	{Pure Math.\ Dept.\\
	University of Waterloo\\
	Waterloo, Ontario\\
	Canada  \  N2L 3G1}
		
\thanks{${}^1$ Research supported in part by NSERC (Canada)}
\thanks{{\ifcase\month\or Jan.\or Feb.\or March\or April\or May\or 
June\or
July\or Aug.\or Sept.\or Oct.\or Nov.\or Dec.\fi\space \number\day,
\number\year}}

\begin{abstract}
We provide a negative solution to a question of M.~Rieffel who asked if the right and left topological stable ranks of a Banach algebra must always agree.   Our example is found amongst a class of nest algebras.   We show that for many other nest algebras, both the left and right topological stable ranks are infinite.   We extend this latter result to Popescu's non-commutative disc algebras and to free semigroup algebras as well.
\end{abstract}

\maketitle
\markboth{\textsc{K.R. Davidson, R. Levene, L.W. Marcoux and H. Radjavi}}{\textsc{On the topological stable rank of non-selfadjoint operator algebras}}


\section{Introduction}
The study of topological stable rank for Banach algebras originated with Rieffel~\cite{Rie1983}.   Motivated by a search for stability results in $C^*$-algebras, he introduced topological stable rank as a non-commutative analogue of the covering dimension for compact spaces.

Given a unital Banach algebra $\cA$, we denote by $Lg_n(\cA)$ (resp.\ $Rg_n(\cA))$ the set of $n$-tuples of elements of $\cA$ which generate $\cA$ as a left ideal (resp.\ as a right ideal).   That is, $Lg_n(\cA) = \{ (a_1, a_2, ..., a_n): a_i \in \cA, 1\le i \le n \mbox{ and there exists } b_1, b_2, ..., b_n \in \cA \mbox{ such that } \sumin b_i a_i = 1 \}$.  The \emph{left} (resp.\ \emph{right}) \emph{topological stable rank} of $\cA$, denoted by $\ltsr(\cA)$ (resp.\ $\rtsr(\cA)$), is the least positive integer $n$ for which $Lg_n(\cA)$ (resp.\ $Rg_n(\cA)$) is dense in $\cA^n$.   When no such integer exists, we set $\ltsr(\cA) = \infty$ (resp.\ $\rtsr(\cA) = \infty$).  If $\ltsr(\cA) = \rtsr(\cA)$, we refer to their common value simply as the \emph{topological stable rank} of $\cA$, written $\tsr(\cA)$.   If $\cA$ is not unital, we define the left (resp.\ the right) topological stable rank of $\cA$ to be that of its unitization.

For $C^*$-algebras, it was shown by Herman and Vaserstein~\cite{HV1984} that topological stable rank coincides with the ring-theoretic notion of stable rank, first introduced by Bass~\cite{Bas1964}.   Consider a ring $\cR$ with identity.   The \emph{left Bass stable rank} of $\cR$, $\lbsr(\cR)$, is the least positive integer $m$ so that for each $(a_1, a_2, ..., a_{m+1}) \in Lg_{m+1}(\cR)$, there exists $(b_1, b_2, ..., b_m) \in \cR^m$ for which $\sum_{i=1}^m (a_i + b_i a_{m+1}) \in Lg_m(\cR)$.   The right Bass stable rank of $\cR$ is analogously defined.   Vaserstein~\cite{Vas1971} (see also Warfield~\cite{War1980}) has shown that $\lbsr(\cR) = \rbsr(\cR)$ for all rings, and hence one normally speaks only of \emph{Bass stable rank} $\bsr(\cR)$.

For general Banach algebras $\cA$ we have the inequality
\[
\bsr(\cA) \le \min( \ltsr(\cA), \rtsr(\cA))\]
(see Corollary~2.4 of~\cite{Rie1983}).    Jones, Marshall and Wolff~\cite{JMW1986} have shown that the disc algebra $\cA(\bbD)$ satisfies $\bsr(\cA(\bbD)) = 1$.  (Recall that the disc algebra $\cA(\bbD)$ consists of those functions which are continuous on the closed unit disc of $\bbC$ and which are analytic on the open unit disc.)   Rieffel~\cite{Rie1983} had shown that $\tsr(\cA(\bbD)) = 2$.  This shows that the inequality above may be strict.

Question~1.5 of Rieffel's paper asks whether or not there exists a Banach algebra $\cA$ for which $\ltsr(\cA) \not =  \rtsr(\cA)$.   It is clear that if such an algebra is to exist, there must be something inherently different between the structure of the left and of the right ideals of $\cA$.   If $\cA$ is a $C^*$-algebra, then the involution provides an anti-isomorphism between left and right ideals, and so one would expect that $\rtsr(\cA)$ should equal $\ltsr(\cA)$ for these algebras.   That this is the case is the conclusion of Proposition~1.6 of~\cite{Rie1983}.

Thus, the search for an algebra $\cA$ of Hilbert space operators for which the left and right topological stable ranks differ takes us into the class of non-selfadjoint algebras.   Two of the best studied such classes are nest algebras and free semigroup algebras.  

We begin the second section of this paper by presenting an example of a nest algebra $\cT(\cN)$ for which $\ltsr (\cT(\cN)) = \infty$ while $\rtsr (\cT(\cN)) = 2$.  The example is found amongst atomic nest algebras, order isomorphic to the natural numbers $\bbN$, all of whose atoms are finite dimensional, with the dimensions of the atoms growing sufficiently rapidly.   We then examine the left and right topological stable ranks of nest algebras in general, and show that in many other cases, the stable ranks agree and are infinite.  It is not yet clear which nest algebras satisfy $\ltsr(\cT(\cN)) = \rtsr(\cT(\cN))$, or indeed, which values of the left (or right) topological stable ranks are attainable. 

In the third section, we deal with the case of  non-commutative disc algebras and of \textsc{wot}-closed free semigroup algebras.  We show that the left and right topological stable ranks of such algebras are always infinite.

\bigskip

First let us prepare the groundwork for what will follow.  We shall need the following two results  due to Rieffel~\cite{Rie1983}.

\begin{thm} \label{Rieffel}
\begin{itemize}
	\item[(a)]
	Let $\cA$ be a Banach algebra and let $J$ be an ideal of $\cA$.    Then $\ltsr(A/J) \leq \ltsr(\cA)$.
	\item[(b)]
	Let $\fH$ be an infinite dimensional, complex Hilbert space.  Then $\tsr(\cB(\fH)) = \infty$.
\end{itemize}	
\end{thm}

A trivial modification of Theorem~\ref{Rieffel}(a) shows that if $\varphi: \cA \to \cB$ is a continuous unital homomorphism with dense range, then $\ltsr(\cB) \le \ltsr(\cA)$.   (See Proposition~4.12 of~\cite{Bad1998} for a version of this for topological algebras.)   Of course there is a corresponding result for right topological stable rank.  We shall also require a slightly more general version of Proposition~1.6 of~\cite{Rie1983}.   Its proof is essentially the same as the proof of that Proposition.  

\begin{lem} \label{involution}
Suppose that $\cA$ is a unital Banach algebra with a continuous involution.   Let $\cB$ be a unital (not necessarily selfadjoint) subalgebra of $\cA$.   Then $\ltsr(\cB) = \rtsr(\cB^*)$.  Hence $\rtsr(\cB) = \ltsr(\cB^*)$.
\end{lem}

\begin{pf}
A simple calculation shows that $(b_1, b_2, ..., b_n) \in Lg_n(\cB)$ if and only if $(b_1^*, b_2^*, ..., b_n^*) \in Rg_n(\cB^*)$.   From this the result easily follows.
\end{pf}

When $\cA$ is a subalgebra of operators on a Hilbert space $\fH$ (or on any other vector space for that matter), we may construct row spaces $\cR_n(\cA)$ and column spaces $\cC_n(\cA)$ of $n$-tuples of elements of $\cA$.   This allows us to view elements of $\cR_n(\cA)$ as operators from $\fH^{(n)}$ to $\fH$, and to view elements of $\cC_n(\cA)$ as operators from $\fH$ to $\fH^{(n)}$.   To say that an element $A = \begin{bmatrix} A_1 & A_2 & \cdots & A_n \end{bmatrix}$ of $\cR_n(\cA)$ lies in $Rg_n(\cA)$ is equivalent to saying that $A$ is right invertible, i.e. that there exists $B = \begin{bmatrix} B_1 & B_2 & \cdots & B_n\end{bmatrix}^t \in \cC_n(\cA)$ such that $A B$ is the identity operator on $\fH$.   That there exists a corresponding statement for $Lg_n(\cA)$ is clear.

Our main tool for determining the topological stable ranks of the algebras considered below is the following observation:

\begin{rem} \label{main_remark}
If an algebra $\cA$ of operators in $\bofh$ contains operators $A_1, ..., A_n$ 
so that $[A_1 \ A_2 \ \cdots \ A_n] \in \cB(\fH^{(n)}, \fH)$ is a semi-Fredholm 
operator of  negative  semi-Fredholm index, then $\rtsr(\cA) \ge n+1$.     
In particular, therefore, if $\cR_n(\cA)$ 
contains a proper isometry, then $\rtsr(\cA) \ge n+1$.
This follows from basic Fredholm theory (see, for eg.~\cite{CPY1974}), as no small perturbation 
$[A_1^\prime \ A_2^\prime \ \cdots \ A_n^\prime]$ of $[A_1 \ A_2 \ \cdots \ A_n] $ 
will be surjective, and thus $\sum_{i=1}^n A_i^\prime B_i \not = I$ 
for any choice of $B_1, B_2, ..., B_n \in \cA$.

The corresponding result for left topological stable rank says that if $\cC_n(\cA)$ contains a proper co-isometry, then $\ltsr(\cA) \ge n+1$.

\end{rem}

The way this observation will be used is as follows:

\begin{prop}  \label{two_isometries}
Suppose that $\cA \subseteq \bofh$ is a Banach algebra of operators and that $\cA$ contains two isometries $U$ and $V$ with mutually orthogonal ranges.   Then $\rtsr(\cA) = \infty$.
\end{prop}

\begin{pf}
Once  $\cA$ contains two such isometries $U$ and $V$, it is clear that for each $n \ge 1$,  $\{ U, VU, V^2 U, ..., V^n U \}$ are $n+1$ isometries in $\cA$ with mutually orthogonal ranges.   Let $Y = [U \ VU \ V^2 U \ \cdots \ V^{n-1}U ] \in \cB(\hilb^{(n)}, \hilb)$; then $Y$ is an isometry and $\ran\, Y$ is orthogonal to $\ran\, V^n U$, so that $Y$ is in fact a proper isometry.  

By Remark~\ref{main_remark}, $\rtsr(\cA) \ge n+1$.  Since $n \ge 1$ was arbitrary, $\rtsr(\cA) = \infty$.
\end{pf}

Of course, if $\cA$ contains two co-isometries with mutually orthogonal initial spaces, then by considering $\cB = \cA^*$, we get $\ltsr(\cA) = \rtsr(\cB) = \infty$.


\section{Nest algebras} \label{nest_algebras}

The first class of algebras we shall examine are \emph{nest algebras}, which are an infinite dimensional generalization of the algebra $\cT_n(\bbC)$ of upper triangular $n \times n$ matrices.   A \emph{nest} $\cN$ on a Hilbert space $\fH$ is a chain of closed subspaces of $\fH$ such that $\{ 0 \}, \fH$ lie in $\cN$, and $\cN$ is closed under the operations of taking arbitrary intersections and closed linear spans of its elements.    At times it is convenient to identify the nest $\cN$ with the collection $\cP(\cN) = \{ P(N): N \in \cN \}$, where - for a subspace $M$ of $\fH$, $P(M)$ denotes the orthogonal projection of $\fH$ onto $M$.    For each $N \in \cN$, we may define the \emph{successor} of $N$ to be $N_+ := \inf \{ M \in \cN: M > N \}$.  If $N_+ \not = N$, then $N_+ \ominus N$ is called an \emph{atom} of $\cN$.   If $\fH$ is spanned by the atoms of $\cN$, we say that $\cN$ is \emph{atomic}.   If $\cN$ admits no atoms, we say that $\cN$ is \emph{continuous}.   Most nests are neither atomic nor continuous.

Given a nest $\cN$, there corresponds to $\cN$ the (\textsc{wot}-closed) \emph{nest algebra}
\[
\cT(\cN) = \{ T \in \cB(\fH): T N \subseteq N \mbox{ for all } N \in \cN \}. \]

A very important example of a nest algebra is the following.   Suppose that $\fH$ is a separable Hilbert space with orthonormal basis $\{e_n \}_{n=1}^\infty$.   Let $N_0 = \{ 0\}$, $N_\infty = \fH$, and for $n \ge 1$, let $N_n = \mathrm{span} \{ e_1, e_2, ..., e_n \}$.   Then $\cN = \{ N_k: 0 \le k \le \infty \}$ is a nest.   
The corresponding nest algebra coincides with the set of all operators in $\cB(\fH)$ whose matrix with respect to this orthonormal basis is upper triangular.   Because of the obvious relation of this nest to the natural numbers, we shall denote this nest algebra by $\cT(\bbN)$.  It is also worth noting that if $\cN$ is a nest on $\fH$, then so is $\cN^\perp := \{ N^\perp: N \in \cN \}$.   In fact, $\cT(\cN^\perp) = \{ T^*: T \in \cT(\cN) \} = \cT(\cN)^*$.   We denote by $\cD(\cN) = \cT(\cN) \cap \cT(\cN)^*$ the \emph{diagonal} of $\cT(\cN)$.  This is a von Neumann algebra.   If $\cN$ is atomic, then it is known that there exists a unique expectation of $\cT(\cN)$ onto $\cD(\cN)$ (see, for eg., Chapter~8 of~\cite{Dav1988}).

\bigskip

The following is the main result of the paper.   It provides an example of a Banach algebra for which the right and left topological stable ranks differ, thereby answering Question~1.5 of~\cite{Rie1983} in the negative.  We thank J.~Orr for simplifying one of the calculations at the end of the proof.
 
\begin{thm}  \label{main_eg}
Let $\cN$ be an atomic nest which is order isomorphic to $\bbN$, 
with finite dimensional atoms $E_k = P(N_k)-P(N_{k-1})$ of 
rank $n_k$ satisfying $n_k \ge 4 \sum_{i<k} n_i$.   
Then 
\[ \ltsr(\cT(\cN)) = \infty \qand \rtsr(\cT(\cN)) =2 .\]
\end{thm}

\begin{pf}
Let $\{ e_{k j} :  1 \le j \le n_k \}$ be an orthonormal basis for the atom $E_k$, $k \ge 1$.   We can construct two co-isometries $U$ and $V$ with mutually orthogonal initial spaces in $\cT(\cN)$ by defining $U^* e_{k j} = e_{2^k 3^j   \ 1}$ and $V^* e_{k j} = e_{5^k 3^j \ 1}$ for all $1 \le j \le n_k$ and $k \ge 1$.
By the remark following Proposition~\ref{two_isometries}, $\ltsr(\cT(\cN)) = \infty$.

It is a consequence of Proposition~3.1 of Rieffel~\cite{Rie1983}, that $\rtsr(\cT(\cN)) \ge 2$.

Let $\Delta$ be the expectation $\Delta(A) = \sum_{k\ge1} E_kAE_k$ of $\cT(\cN)$
onto the diagonal $\cD(\cN)$, which is a finite von Neumann algebra.
Every element $D \in \cD(\cN)$ factors as $D=UP$
where $P$ is positive and $U$ is unitary.
Thus for any $\eps>0$, $D'=U(P+\eps I)$ is an $\eps$-perturbation
which is invertible with inverse bounded by $\eps^{-1}$.

Let $A$ and $B$ belong to $\cT(\cN)$, and let $\eps>0$ be given.
By the previous paragraph, there are $\eps/2$-perturbations $A', B'$
of $A$ and $B$ so that 
$A'= D_a + A'_0$ and $B'=D_b + B'_0$ where $A'_0, B'_0$ lie
in the ideal $\cT_0(\cN)$ of strictly upper triangular operators
and $D_a, D_b$ are invertible elements of $\cD(\cN)$ 
with inverses bounded by $2\eps^{-1}$.
Let
\[ A_1 = A' D_a^{-1} = I + A_0 \qand B_1 = B'D_b^{-1} = I + B_0 ,\]
where $A_0 = A'_0D_a^{-1}$ and $B_0 = B'_0D_b^{-1}$.

Now $A_0 = \sum_{k\ge2} A_0E_k$ and 
\[
 \rank(A_0 E_k) = \rank(P(N_{k-1}) A_0 E_k) 
 \le \rank(P(N_{k-1})) = \sum_{i<k}n_i \le \frac{n_k}4 .
\]
The same estimate holds for $B_0$.  
Therefore we may select projections $P_k \le E_k$
with $\rank P_k \le n_k/2$
so that $A_0E_k=A_0P_k$ and $B_0E_k=B_0P_k$. 

Let $U_k = P_kU_k(E_k-P_k)$ be a partial isometry 
with range $P_k\hilb$.
Define operators $U = \sum_{k\ge1}U_k$ and $P=\sum_{k\ge1}P_k$;
so $P^\perp = \sum_{k\ge1}(E_k - P_k)$ and $UU^* = P$.
Take any positive number 
\[ 0 < \delta < \frac{\eps}2  \|D_a\|^{-1} .\]
Consider $A'' = A' + \delta UD_a$.
Then 
\[ \|A-A''\| < \|A-A'\| + \delta\|D_a\| < \eps ,\]
and from above, $\| B - B' \| < \eps/2$.
We will show that $\big[ A''\ B' \big]$ is right invertible.


Observe that $A' D_a^{-1} P^\perp = B' D_b^{-1} P^\perp = P^\perp$.      Thus
\begin{align*}
A'' D_a^{-1} P^\perp U^* + B' D_b^{-1} P^\perp (I - U^*) &=
P^\perp U^* + \delta U D_a D_a^{-1} P^\perp U^* + P^\perp(I-U^*)\\ &=
\delta P + P^\perp.
\end{align*}
It is clear that this is right invertible (by $\begin{bmatrix} \delta^{-1} P \\ P^\perp \end{bmatrix}$), whence $\big[ A'' \ B' \big]$ is also right invertible.

It follows that $Rg_2(\cT(\cN))$ is dense in $\cR_2(\cT(\cN))$;
that is, $\rtsr(\cT(\cN) ) = 2$.
\end{pf}


Let us next turn our attention to general nest algebras. 
We can show in a large number of cases, the left and right topological stable ranks of a nest algebra agree, and that they are infinite.  For the remainder of this article, we shall restrict our attention to \emph{complex, infinite dimensional, separable Hilbert spaces}.

\begin{prop} \label{decreasing_sequence}
Let $\cN$ be a nest on a Hilbert space $\fH$, and suppose that $\cN$ contains a strictly decreasing sequence $\{ N_k \}_{k=0}^\infty$.   Then $\rtsr(\cT(\cN)) = \infty$.
\end{prop}

\begin{pf}
Let $N_\infty = \cap_{k \ge 0} N_k \in \cN$.  If $\fK := N_0 \ominus N_\infty$, then $\cM = \{ N \cap \fK: N \in \cN \}$ is a nest, and the compression map
\[
\begin{array}{rccc}
\Gamma: & \cT(\cN) & \to & \cT(\cM) \\
 & T & \mapsto & T_{|{\fK}}
\end{array} \]
is a contractive, surjective homomorphism of $\cT(\cN)$ onto $\cT(\cM)$.   By Theorem~\ref{Rieffel}, it suffices to prove that $\rtsr(\cT(\cM)) = \infty$.

If $M_k := N_k \ominus N_\infty$, then $M_k \in \cM$ for all $k \ge 1$, and $M_0 > M_1 > M_2 > \cdots.$   Let $A_k = M_{k-1} \ominus M_k$, $k \ge 1$, and choose an orthonormal basis $\{ e_{k j}: 1 \le j < n_k \}$ for $A_k$, where $2 \le n_k \le \infty$.   Observe that $\cup_{k \ge 1}  \{ e_{k j}: 1 \le j < n_k \}$ is then an orthonormal basis for $\fK$.   

We then define two isometries $U, V \in \cT(\cM)$ via:
\[
U e_{k j} = e_{2^j 3^k \ 1}, \ \ \ \ \ \ \ \ \ \ \ \ \ \ \ \ \ \ \ \ V e_{k j} = e_{5^j 7^k \ 1}, \]
for all $1 \le j < n_k$, $1 \le k < \infty$.   Clearly  $U$ and $V$ have mutually orthogonal ranges.   By Proposition~\ref{two_isometries}, $\rtsr(\cT(\cM)) = \infty$, which -- as we have seen -- ensures that $\rtsr(\cT(\cN)) = \infty$.
\end{pf}


\begin{cor} \label{increasing_sequence}
Let $\cN$ be a nest on a Hilbert space $\fH$, and suppose that $\cN$ contains a strictly increasing sequence $\{ N_k \}_{k=0}^\infty$.   Then $\ltsr(\cT(\cN)) = \infty$.
\end{cor}

 
In the following Theorem, we refer to the \emph{dual} of an ordinal.  If $(\beta, \le)$ is an ordinal, the dual of $\beta$ is the totally ordered set $(\beta^*, \le_*)$ where $\beta^* = \beta$ and $x \le_* y$ if and only if $y \le x$.

\bigskip

\begin{thm} \label{general}
Let $\cN$ be a nest on a Hilbert space $\fH$.   If $\cN$ satisfies any one of the following three properties, then $\ltsr(\cT(\cN)) = \rtsr(\cT(\cN)) = \infty$.
\begin{enumerate}
	\item[(a)]
	$\cN$ has an infinite dimensional atom.
	\item[(b)]
	$\cN$ is uncountable.
	\item[(c)]
	$\cN$ is countable, but is not order isomorphic to an ordinal or its dual.
\end{enumerate}
\end{thm}

\begin{pf}
(a) Choose $N \in \cN$ so that $\dim (N_+ \ominus N) = \infty$, and set $E = N_+ \ominus N$.
Then $E$ is a semi-invariant subspace for $\cT(\cN)$, and the map 
$\gamma: \cT(\cN) \to \cB(E)$ defined by $X \mapsto P(E)X|_E$ is a surjective homomorphism.    Now $\mathrm{tsr}(\cB(E)) = \infty$, by Theorem~\ref{Rieffel}(b).   
Furthermore, by Theorem~\ref{Rieffel}(a), since $\cB(E)$ is a homomorphic image of 
$\cT(\cN)$, $\mathrm{ltsr}(\cT(\cN)) \geq \mathrm{ltsr}(\cB(E)) = \infty$,
and similarly $\mathrm{rtsr}(\cT(\cN)) \geq \mathrm{rtsr}(\cB(E)) = \infty$, 
completing the proof.
\medbreak

(b,c) In each of these cases, the conditions on $\cN$ guarantee the existence of both a strictly increasing sequence $\{ N_k\}_{k=1}^\infty$ and a strictly decreasing sequence 
$\{ M_k \}_{k=1}^\infty$ of subspaces in $\cN$.  The result now follows immediately from 
Proposition~\ref{decreasing_sequence} and Corollary~\ref{increasing_sequence}.
\end{pf}


\begin{cor} \label{maximum_tsr_is_infinite}
Let $\cN$ be a  nest.   Then $\max ( \ltsr(\cT(\cN), \rtsr(\cT(\cN)) = \infty$.
\end{cor}

\begin{pf}
Taking into account the above results, the only case left to consider is that where $\cN$ is a countably infinite nest, order isomorphic to an ordinal or the dual of an ordinal.  As such, $\cN$ either contains an interval $[N_1, N_\infty)$ which is order isomorphic to $\bbN$, or an interval $(N_{-\infty}, N_{-1}]$ which is order isomorphic to $-\bbN$.  From Proposition~\ref{decreasing_sequence} and Corollary~\ref{increasing_sequence} we deduce that $\max ( \ltsr(\cT(\cN), \rtsr(\cT(\cN)) = \infty$.
\end{pf}


We have thus reduced the problem of determining the topological stable ranks of nest algebras to the problem of determining the right topological stable rank of a countable, atomic nest $\cN$, order isomorphic to an ordinal, all of whose atoms are finite dimensional.   Theorem~\ref{main_eg} shows that in this case it is possible to have $\rtsr(\cT(\cN)) = 2$.      We shall see below that this example may be extended to a more general class of nest algebras whose nests are totally ordered like $\omega$ (the first infinite ordinal), and for which there is an arithmetically increasing  sequence of atoms whose ranks grow geometrically fast (see Theorem~\ref{extension}).  Having said this, the exact nature of the nests for which the right topological stable rank is finite is not completely understood.   We begin by establishing a couple of conditions on a nest  $\cN$ which will guarantee that the right topological stable rank  of $\cT(\cN)$ is infinite.

The proofs of the results depend upon the existence of certain surjective homomorphisms of nest algebras established in~\cite{DHO1995}.   Since they play such a key role, we briefly recall the construction of these homomorphisms as outlined in that paper.


\subsection{The Davidson-Harrison-Orr Construction} \label{construction}
Let $\Omega$ be an interval of $\bbZ$, and suppose that $\Omega = \cup_{n=1}^\infty \Omega_n$, where $\Omega_n \subseteq \Omega_{n+1}$, $n \ge 1$ are subintervals of $\Omega$.   Suppose also that $\cM$ is a nest, order isomorphic to $\Omega$ via an order isomorphism $\lambda$.   Let $E_n$ denote the subinterval of $\cM$ corresponding via $\lambda$ to the interval $\Omega_n$, $n \ge 1$.

Consider next a nest $\cN$ containing countably many subintervals $F_n$ acting on pairwise orthogonal subspaces such that $\cT(\cN)|_{F_n \fH}$ is unitarily equivalent to $\cT(\cM)|_{E_n \fH}$ via a unitary conjugation $Ad_{U_n}: \cT(\cN)|_{F_n \fH} \to \cT(\cM)|_{E_n \fH}$.   If $\alpha_n: \cT(\cN) \to \cT(\cN)|_{F_n \fH}$  are the natural compression maps and $\beta_n: \cT(\cM)|_{E_n \fH} \to \cT(\cM)$ are the inclusion maps, $ n \ge 1$, let $\varphi_n: \cT(\cN) \to \cT(\cM)$ be the maps $\varphi_n = \beta_n \circ Ad_{U_n} \circ \alpha_n$, $n \ge 1$, so that $\varphi_n$ is a homomorphism for all $n$.

Letting $\cU$ be a free ultrafilter on $\bbN$, we have that
\[
\begin{array}{rccc}
\varphi: & \cT(\cN) & \to & \cT(\cM) \\
\ & T & \mapsto & \textsc{wot}-\lim_{n \in \cU} \varphi_n(T)
\end{array} \]
defines a continuous epimorphism of $\cT(\cN)$ onto $\cT(\cM)$ (\cite{DHO1995}, Corollary~5.3 and Theorem~6.8).

For example, suppose that $\cM$ is the maximal atomic nest, ordered like $\omega^*$, so that $\cT(\cM) \simeq \cT(\bbN)^*$.     Set $\Omega = -\bbN$, $\Omega_n = \{ -2^n, -2^n +1, ..., -3, -2, -1 \}$, and let $E_n$ denote the corresponding subinterval of $\cM$.   Thus $\cT(\cM)|_{E_n \fH} \simeq \cT_{2^n}(\bbC)$, the upper triangular $2^n \times 2^n$ matrices over $\bbC$.   Choose integers $r_1 < r_2 < r_3 < \cdots$  such that $r_n - r_{n-1} > 2^n$.   If $F_n = \mathrm{span} \{ e_{r_n +1}, e_{r_n +2}, ..., e_{r_n + 2^n} \}$, then the $F_n$'s are pairwise orthogonal and $\cT(\bbN)|{F_n \fH} \simeq \cT_{2^n}(\bbC)$ as well, and so we can find a unitary matrix $U_n: F_n \fH \to E_n \fH$ such that $\cT(\bbN)|_{F_n\fH} = U_n^* (\cT(\cM)|_{E_n \fH}) U_n$.     With $\cU$ a free ultrafilter on $\bbN$, 
\[
\varphi(T) = \textsc{wot}-lim_{n \in \cU} \varphi_n(T) \]
implements a continuous epimorphism of $\cT(\bbN)$ onto $\cT(\cM) \simeq \cT(\bbN)^*$.


\begin{cor} 
$\ltsr(\cT(\bbN)) = \rtsr(\cT(\bbN)) = \infty$.
\end{cor}

\begin{pf}
By Corollary~\ref{increasing_sequence}, $\ltsr(\cT(\bbN)) = \infty$.  Let $\varphi: \cT(\bbN) \to \cT(\bbN)^*$ be the epimorphism described in the Section~\ref{construction}.    By Theorem~\ref{Rieffel} and Lemma~\ref{involution}, $\rtsr(\cT(\bbN)) \ge \rtsr(\cT(\bbN)^*) = \ltsr(\cT(\bbN)) = \infty$.
\end{pf}


\begin{rem} \label{remark}
More generally, suppose that $\cN$ is a countable nest, order isomorphic to an ordinal, and that $\cN$ contains intervals of length $n_1 < n_2 < n_3 < \cdots$, such that the interval with length $n_k$ has consecutive atoms of  size $(d_{n_k}, d_{n_{k-1}}, ..., d_1)$.   Without loss of generality, we may assume that the subspaces upon which these intervals act are mutually orthogonal.  The above construction can be used to produce an epimorphism of $\cT(\cN)$ onto $\cT(\cM)$, where $\cM$ is a nest of order type $\omega^*$ (and whose atoms have dimensions $(..., d_{4}, d_{3}, d_{2}, d_{1})$).   By Proposition~\ref{decreasing_sequence}, $\rtsr(\cT(\cM)) = \infty$, and thus by Theorem~\ref{Rieffel}, $\rtsr(\cT(\cN)) = \infty$ as well.
\end{rem}


\begin{eg}
Let $\cN$ be the nest order isomorphic to $\omega$, whose atoms $(A_n)_{n=1}^\infty$ have dimensions $1, 2, 1, 3, 2, 1, 4, 3, 2, 1, ...., n, n-1, n-2, ..., 3, 2, 1, n+1, n, n-1, ...$.  Then $\ltsr(\cT(\cN)) = \rtsr(\cT(\cN)) =\infty$.
\end{eg}

\smallskip

The next result is an immediate consequence of the Remark~\ref{remark}.


\begin{cor} \label{bounded_size}
Let $\cN$ be a countable nest, order isomorphic to an ordinal.   Suppose that $\cN$ contains intervals $E_j$ of length $n_j$, where $n_j < n_{j+1}$ for all $j$, and such that $\max \{ \dim\, A: A \in E_j \mbox{ an atom} \} < K$ for some constant $K \ge 1$ independent of $j$.  Then $\rtsr(\cT(\cN)) = \infty$.
\end{cor}

We remark that in Remark~\ref{remark} and in Corollary~\ref{bounded_size}, the assumption that $\cN$ be countable and order isomorphic to an ordinal is stronger than what is needed to obtain an epimorphism of $\cT(\cN)$ onto a nest algebra $\cT(\cM)$ with right topological stable rank equal to $\infty$.  On the other hand, it simplifies the exposition, and the right topological stable rank of $\cT(\cN)$ in all other cases has been dealt with already.

\vskip 0.5 cm

Theorem~\ref{main_eg} shows that if $\cN$ is a nest, ordered like the natural numbers, whose atoms grow geometrically fast in dimension, then $\rtsr(\cT(\cN)) = 2$.   The conditions on the rate of growth of the dimensions of the atoms can be somewhat relaxed.  The following observation will prove useful.

\bigskip

Let $B = \begin{bmatrix} B_1 & B_2 \\ 0 & B_4 \end{bmatrix} \in \cB(\fH_1 \oplus \fH_2)$ be an operator where $B_1, B_4$ are invertible.   Then $B$ is invertible with $B^{-1} = \begin{bmatrix} B_{1}^{-1} & -B_{1}^{-1} B_2 B_4^{-1} \\ 0 & B_4^{-1} \end{bmatrix}$.   Thus if there exist a constant $H  > 0$ so that $\norm B_1^{-1} \norm \le H$, $\norm B_4^{-1} \norm \le H$, then $\norm B^{-1} \norm \le 2 H + H^2 \norm B \norm$.

Using induction, it is not hard to see that if 
\[
 A = \begin{bmatrix}
  A_{1 1} & A_{1 2} & \ldots & A_{1 n} \\  
  \ & A_{2 2} & \ldots & A_{2 n} \\ 
  & & \ddots & \\ & & & A_{n n} 
 \end{bmatrix}
\]
is an operator in $\cB(\oplus_{k=1}^n \fH_k)$ and if each $A_{kk}$ is invertible with $\norm A_{kk}^{-1} \norm \le H$ for some $H > 0$, then $\norm A^{-1} \norm \le L$ for some constant $L$ that depends only upon $H, \ n$ and $\norm A \norm$.


\begin{thm} \label{extension}
Suppose that $\cN$ is a nest ordered like $\omega$, all of whose atoms are finite dimensional.   Let $A_n$, $n \ge 1$, denote the atoms of $\cN$, and let $r_n = \dim\, A_n$ for $n\ge1$. 
Set $R(k) = \max_{1 \le i \le  k} r_i$ for $k \ge 1$.
Suppose that there exists a $\gamma >0$ and an integer $J > 0$  such that 
\[ R((k+1)J) \ge (1 + \gamma) R(kJ) \quad\text{for all}\quad k \ge 1. \]
Then 
\[ \ltsr(\cT(\cN)) = \infty \quad\text{and}\quad \rtsr(\cT(\cN)) = 2. \]
\end{thm}

\begin{pf}
By Corollary~\ref{increasing_sequence}, $\ltsr(\cT(\cN)) = \infty$.   

Suppose that $\gamma > 0$ and $J \ge 1$ are chosen as in the statement of the Theorem, and that 
\begin{align*}
R((k+1)J) 
	&= \max\, \{ r_k: 1 \le i \le (k+1)J \} \\
	&= \max \{ r_i: k J \le i \le (k+1) J \}  
	\ge (1 + \gamma) R(kJ) .
\end{align*}

Choose an integer $p \ge 1$ so that $\displaystyle \frac{(1 + \gamma)^p}{p} \ge 5 J$.  
For $k \ge 1$, set 
\[ F_k = \sum_{(k-1) p J < i \le k p J } P(A_i)   .\]
In essence, we are grouping together blocks of length $p J$ of $\cN$ into a single ``megablock".  The compression of $\cT(\cN)$ to any such ``megablock" is an upper triangular $p J \times p J$ operator matrix whose entries are finite dimensional matrices.

Observe that if $M_k := \rank F_k$, then $M_k \le p J R(kpJ)$ 
(since the maximum rank among the atoms of $F_k$ is $R(kpJ)$) and that 
\begin{align*}
 M_{k+1} &= \rank F_{k+1} \ge R((k+1)pJ) \\ 
 &\ge (1 + \gamma)^p R(kpJ) \ge 5 p J R(k p J) \ge 5 M_k . 
\end{align*}
Hence $M_k \ge 4 \sum_{i < k} M_i$ for each $k \ge 2$.	

The remainder of the proof will be an adaptation of the proof of Theorem~\ref{main_eg};  the main difference being that we will work with blocks of length $p J$ of $\cN$.

Let $A, B \in \cT(\cN)$ be given, and let $\varepsilon > 0$.   Let $\Delta(A) = \sum_{k \ge 1} E_k A E_k$ be the expectation of $A$ onto the diagonal $\cD(\cN)$ of $\cT(\cN)$.   As noted in the proof of Theorem~\ref{main_eg}, $\cD(\cN)$ is a finite von Neumann algebra and so $\Delta(A) = U P$ for some unitary $U$ and positive operator $P$ lying in $\cD(\cN)$.  But then $\Delta(A)^\prime = U (P + \frac{\varepsilon}{2} I)$ is an $\varepsilon/2$-perturbation of $\Delta(A)$ which is invertible with inverse bounded above by $H = \frac{2}{\varepsilon}$.   

Let $A^\prime = \Delta(A)^\prime + (A - \Delta(A))$.   Note that the compression of $A^\prime$ to $F_k \fH$ is a $p J \times p J$ block-upper triangular matrix whose diagonal entries are all invertible with inverses bounded above by $H$.   By the comments preceding this Theorem, $F_k A^\prime F_k$ is invertible with $\norm (F_k A^\prime F_k)^{-1} \norm \le L_A$, where $L_A$ is a constant depending only upon $\varepsilon,  p J, $ and $\norm A \norm$.   A similar construction applied to $B$ yields an operator $B^\prime$ such that $F_k B^\prime F_k$ is invertible with $\norm (F_k B^\prime F_k)^{-1} \norm \le L_B$ for all $k \ge 1$, where $L_B$ is a constant depending only upon $\varepsilon, p J$ and $\norm B \norm$.

Thus we can write $A^\prime = D_a + A_0^\prime$, $B^\prime = D_b + B_0^\prime$, where $D_a = \sum_{k \ge 1} F_k A^\prime F_k$, $D_b = \sum_{k \ge 1} F_k B^\prime F_k$ are invertible elements of $\cD = \sum_{k \ge 1} (F_k \cT(\cN) F_k)$, $\norm D_a^{-1} \norm  \le L_A$, $\norm D_b^{-1} \norm \le L_B$ and $A^\prime_0 := A^\prime - D_a$, $B^\prime_0 := B^\prime -D_b$ lie in the ideal 
\[ \cT^\prime_0(\cN) = \{ T \in \cT(\cN): \sum_{k \ge 1} F_k T F_k = 0 \} .\]

At this point, we can apply the second half of the proof of Theorem~\ref{main_eg}, with the $E_k$'s of that proof replaced with $F_k$, $k \ge 1$, to conclude that $\rtsr(\cT(\cN)) = 2$.
\end{pf}

One very interesting consequence of Theorems~\ref{main_eg} and~\ref{extension} is that they allow us to resolve (in certain cases) a question of Davidson, Harrison and Orr (see~\cite{DHO1995}, Section~8) regarding epimorphisms of nest algebras onto $\cB(\fH)$.   

\begin{prop}
Let $\cN$ be a nest  of the type described in Theorem~\ref{extension}.   Let $\cA$ be an operator algebra with $\rtsr(\cA) = \infty$.    Then there is no epimorphism of $\cT(\cN)$ onto $\cA$.   In particular, this holds if $\cA$ is any one of the following:
\begin{enumerate}
	\item[(a)]
	$\cB(\fH)$;
	\item[(b)]
	$\cT(\cV)$, where $\cV$ is an uncountable nest; or 
	\item[(c)]
	$\cT(\cM)$, where $\cM$ is a countable nest which is not isomorphic to an ordinal.
\end{enumerate}	
\end{prop}

\begin{pf}
Observe that $\rtsr(\cB(\fH)) = \rtsr(\cT(\cV)) = \rtsr(\cT(\cM)) = \infty$ by Theorem~\ref{general}.
If such an epimorphism were to exist, then by Theorem~\ref{Rieffel}, it would follow that $\rtsr(\cT(\cN)) = \infty$, which is a contradiction. 
\end{pf}

We finish this section by mentioning a few unresolved questions dealing with the stable rank of nest algebras.

There are still a number of nests for which we have been unable to determine the left and right topological stable ranks.   When the nest is ordered like $\omega$, it is clear that the value of the right topological stable rank of the corresponding nest algebra depends upon how fast the atoms grow.  If the atoms of $\cN$ are bounded in dimension, then $\rtsr(\cT(\cN)) = \infty$.   If the dimensions of the atoms grow at an exponential rate, then the right topological stable rank is $2$.  What happens when the rate of growth lies between these two extremes?  A key case which we have been unable to resolve and which would very likely shed light upon the general problem is the following:

\smallskip

\noindent{\textbf{Question 1.}}  Suppose that $\cN$ is a nest, ordered like $\omega$, whose atoms $(A_n)_{n=1}^\infty$ satisfy $\dim \, A_n = n$, $n \ge 1$.  What is $\rtsr(\cT(\cN))$?   

\smallskip

We note that by Corollary~\ref{increasing_sequence}, $\ltsr(\cT(\cN)) =\infty$.  

\smallskip

Observe that in all of our examples, $\rtsr(\cT(\cN)) \in \{2, \infty\}$.   
\smallskip

\noindent{\textbf{Question 2.}} Does there exist a  countable nest $\cN$, order isomorphic to an ordinal (in particular - order isomorphic to $\omega$), all of whose atoms are finite dimensional, for which the value of $\rtsr(\cT(\cN))$ is other than $2$ or $\infty$?   

\bigskip

The above analysis  suggests that it is not the exact dimensions of the atoms which is significant, but rather the rate at which these dimensions grow.   If $\rtsr(\cT(\cN)) = m$ for some $3 \le m < \infty$, then by a straightforward adapation of Proposition~6.1 of~\cite{Rie1983} to general Banach algebras, $\rtsr(\cT(\cN)) \otimes \bbM_n  = \lceil  (\rtsr(\cT(\cN)) -1)/m \rceil +1$, and hence $\rtsr(\cT(\cN)) \otimes \bbM_n = 2$ for sufficiently large values of $n$.  (Here $\lceil k \rceil$ denotes the least integer greater than or equal to $k$.)  But $\cT(\cN) \otimes \bbM_n \simeq \cT(\cM)$, where $\cM$ is a nest, order isomorphic to $\cN$, whose atoms have dimension $n$ times the dimension of the corresponding atoms of $\cN$.   As such, the rate of growth of the atoms of $\cM$ is identical to that of $\cN$.   We suspect that this should imply that $\rtsr(\cT(\cN)) = 2$, but we have not been able to prove this.

If $\cA$ is any unital Banach algebra with $\rtsr(\cA) = \infty$, then it follows from the previous paragraph that $\rtsr(\cA \otimes \bbM_n) = \infty$ for all $n \ge 1$.  
We obtain the following result for Banach algebras which was established
for C*-algebras by Rieffel~\cite[Theorem 6.4]{Rie1983}.  We shall first fix a basis $\{e_n\}_{n=1}$ for $\fH$, and denote by $E_{i j}$ the matrix unit $e_i e_j^* \in \cK(\fH)$.  
If $\cA$ is a unital Banach algebra, 
consider any Banach algebra cross norm on $\cA \otimes \cK(\fH)$
for which $\cA$ is imbedded isometrically (but not unitally)
as a corner $\cA \otimes E_{11}$, each matrix algebra
$\cA \otimes \bbM_n$ is identified with 
$\big(\sum_{i=1}^n E_{ii} \big) \cA \otimes \cK(\fH) \big(\sum_{j=1}^n E_{jj} \big)$,
and  the union of these matrix algebras is norm dense in 
$\cA \otimes \cK(\fH)$.

\begin{prop}
Let $\cA$ be a Banach algebra with identity.
Then 
\[ \ltsr(\cA \otimes \cK(\fH))=\rtsr(\cA \otimes \cK(\fH)) \in \{1,2\} ,\]
and it equals $1$ if and only if $\tsr(\cA) = 1$.
\end{prop}

\begin{pf}
The argument that the (left or right) topological stable rank is
at most 2 is done by Rieffel~\cite{Rie1983}.  He also shows that $\ltsr(\cA)=1$
and $\rtsr(\cA)=1$ and both equivalent to the density of the
invertible elements.  If the invertibles are dense in $\cA$, Rieffel
shows that they are also dense in $\cA \otimes \bbM_n$ for all
$n \ge 1$.  From this, it is easy to see that the invertibles are
dense in the unitization $(\cA \otimes \cK(\fH))^\sim$.
To complete the proof, it suffices to show that if the 
invertibles are dense in $(\cA \otimes \cK(\fH))^\sim$, 
then they are also dense in $\cA$.

Fix $A \in \cA$ with $\|A\| \le 1/2$.
Let $\cC$ denote the circle centred at $0$ of radius $3/4$.
Define $M = \sup\{ \|(zI-A)^{-1}\| : z \in \cC \} \ge 4$.
Then $A' = A \otimes E_{11} + I \otimes E_{11}^\perp$
belongs to $(\cA \otimes \cK(\fH))^\sim$.
For any $0 < \eps < (6\pi M^2)^{-1} < 1/4$, 
choose $B \in (\cA \otimes \cK(\fH))^\sim$
so that $\|A'-B\| < \eps$.

The spectrum of $A'$ is $\sigma(A') = \sigma(A) \dot\cup \{1\}$.
By \cite[Theorem 1.1]{Her1989}, $\sigma(B)$ is disjoint from $\cC$.
By the Riesz functional calculus, there is an idempotent
\[ P = \int_\cC (zI-B)^{-1} \,dz \]
which commutes with $B$.
This idempotent is close to $E := I \otimes E_{11}$ because of the following
estimates.  For $z \in \cC$,
\begin{align*}
 \| (zI-B)^{-1} \| &= \big\| \big( (zI-A') - (B-A') \big)^{-1} \big\| \\ &= 
 \big\| (zI - A')^{-1} \sum_{n\ge0} \big((B-A')(zI - A')^{-1} \big)^n \big\| \\
 &\le \frac M {1-M\eps} < 2M.
\end{align*}
Therefore
\begin{align*}
 \|P-E\| &= \Big\| \int_\cC   (zI-B)^{-1} - (zI - A')^{-1} \,dz \Big\| \\
 &\le 2 \pi \frac34 \sup_{z \in \cC} \| (zI-B)^{-1} (A'-B) (zI-A')^{-1} \| \\
 &\le \frac {3\pi} 2 2M \eps M = 3\pi M^2 \eps =: \eps' < \frac12 .
\end{align*}

Now a standard argument shows that $S = PE + (I-P)E^\perp$
is an invertible element of  $(\cA \otimes \cK(\fH))^\sim$ such that
$SE = PS$ and 
\[ \|S-I\| = \| (P-E)(E-E^\perp)\| = \|P-E\| \le \eps'.\]
Thus $B' = S^{-1}BS$ is close to $B$ and has the form 
$B' = B_1 \otimes E_{11} + E_{11}^\perp B_2 E_{11}^\perp$.
Indeed,
\begin{align*}
 \|B'-B\|& \le \|S^{-1}\| \, \|(S-I)B - B(S-I)\| 
 \\&\le \frac1{1-\eps'} 2 \|B\| \eps' 
 \le \frac{1+2\eps}{1-\eps'} =: \eps''. 
\end{align*}
Thus we obtain that $\|A-B_1\| < \eps+\eps''$ and $B_1$ is
invertible in $\cA$.  Since $\eps''$ tends to 0 as $\eps$ does,
we conclude that the invertibles are dense in $\cA$.
\end{pf}

\bigskip

Another interesting and open problem concerns the Bass stable rank of nest algebras.   For the nests of Theorem~\ref{main_eg} or more generally for those of Theorem~\ref{extension}, it follows from the inequality mentioned in the introduction that $\bsr(\cT(\cN)) \le \min(\ltsr(\cT(\cN)), \rtsr(\cT(\cN)) = 2$.   Nevertheless, an explicit calculation of $\bsr(\cT(\cN))$ for this or indeed for any nest algebra seems to be an rather difficult problem.  
\smallskip

\noindent{\textbf{Question 3.}} Find $\bsr(\cT(\bbN))$, or indeed $\bsr(\cT(\cN)))$ of any nest algebra.




\section{Non-commutative operator algebras generated by isometries}


Let us now consider operator algebras generated by free semigroups of isometries.   The theory here divides along two lines;  the norm-closed version, often referred to as \emph{non-commutative disc algebras}, and the \textsc{wot}-closed versions, known simply as \emph{free semigroup algebras}.   The latter algebras include the \emph{non-commutative Toeplitz algebras}, to be described below.

Let $n \ge 1$.   The non-commutative disc algebra $\fA_n$, introduced by Popescu~\cite{Pop1991, Pop1996}, is (completely isometrically isomorphic to) the norm-closed subalgebra of $\cB(\fH)$ generated by the identity operator $I$ and $n$ isometries $S_1, S_2, ..., S_n$ with pairwise orthogonal ranges.   It is shown in~\cite{Pop1996} that the complete isometric isomorphism class of $\fA_n$ is independent of the choice of the isometries, and that $\fA_n$ is completely isometrically isomorphic to $\fA_m$ if and only if $m = n$.      Note that for each $1 \le j \le n$, $S_j^* S_j = I \ge \sumin S_i S_i^*$, and when $\sumin S_i S_i^* = I$, the $C^*$-algebra generated by $\{ S_1, S_2, ..., S_n \}$ is the Cuntz algebra $\cO_n$.  When $\sumin S_i S_i^* < I$, the $C^*$-algebra generated by $\{ S_1, S_2, ..., S_n \}$ is the Cuntz-Toeplitz algebra $\cE_n$.

Given isometries $S_1, S_2, ..., S_n$ with pairwise orthogonal ranges as above, the \textsc{wot}-closure $\fS_n$ of the corresponding disc algebra $\fA_n$ is known as a \emph{free semigroup algebra}.    These were first described in~\cite{DP1999}.     Of particular importance is the following example.    Let $\bbF_n^+$ denote the free semigroup on $n$ generators $\{ 1, 2, ..., n \}$.   Consider the Hilbert space $\fK_n = \ell^2(\bbF_n^+)$ with orthonormal basis $\{ \xi_w: w \in \bbF_n^+ \}$.  For each word $v \in \bbF_n^+$, we may define an isometry $L_v \in \cB(\fK_n)$ by setting $L_v \xi_w = \xi_{v w}$ (and extending by linearity and continuity to all of $\fK_n$).   The identity operator is $L_\varnothing$.  Then $L_1, L_2, ..., L_n$ are $n$-isometries with orthogonal ranges, and the \textsc{wot}-closed algebra $\fL_n$ generated by $I, L_1, L_2, ..., L_n$ is called the non-commutative Toeplitz algebra.

A theorem of Davidson, Katsoulis, and Pitts~\cite{DKP2001} shows that if $\fS_n$ is a free semigroup algebra, then there exists a projection $P \in \fS_n$ such that $\fS = \fM P \oplus P^\perp \fM P^\perp$, where $\fM$ is the von Neumann algebra generated by $\fS_n$, and $\fS P^\perp = P^\perp \fS P^\perp$ is completely isometrically isomorphic to $\fL_n$.

\bigskip


\begin{thm} \label{non_comm}
Let $n \ge 2$.   
\begin{enumerate}
	\item[(a)]
	If $\fA_n$  is the non-commutative disc algebra on $n$-generators, then $\tsr(\fA_n) = \infty$.
	\item[(b)]
	If  $\fS_n$   is a free semigroup algebra on $n$-generators, then $\tsr(\fS_n) = \infty$.
\end{enumerate}	
\end{thm}

\begin{pf}
First observe that both $\fA_n$ and $\fS_n$ are generated by $n \ge 2$ isometries with mutually orthogonal ranges.   By Proposition~\ref{two_isometries}, $\rtsr(\fA_n) = \rtsr(\fS_n) = \infty$.

We now consider the left topological stable rank of these two algebras.

Let $V_1, V_2, ..., V_n \in \bofh$ be isometries  with mutually orthogonal ranges.  
Let $A_i = \frac 1 n V_i^*$, $1 \le i \le n$.   
Then $\sum_{i=1}^n A_i^* A_i  = \frac 1 n I$
is a strict contraction. 
By Proposition~2 of~\cite{Bun1984}, 
there exists a Hilbert space $\cK$ containing $\cH$ and pure  isometries 
$\{ W_i \}_{i=1}^n \subseteq \cB(\cK)$  with pairwise orthogonal ranges
so that $\hilb ^\perp  \in \mathrm{Lat}\, W_i$ and 
$ P_{\hilb}  W_i|_{\hilb} = A_i$, $1 \le i \le n$.   

\medbreak
(a) The norm-closed algebra $\cB_n \subseteq \cB(\cK)$ generated by 
$\{ I, W_1, W_2, ..., W_n \}$ satisfies $\cB_n \simeq \fA_n$.   
The compression map
\[
 \begin{array}{rccc}
 \gamma: & \cB_n & \to & \bofh \\
 & X & \mapsto & P_\hilb X|_{\hilb}
 \end{array}	
\]
is a (completely contractive) homomorphism, 
as $\hilb\textcolor{blue}{^\perp} \in \mathrm{Lat}\, W_i$ for all $i$.   
Thus $\rtsr(\cB_n) \ge \rtsr(\overline{\gamma(\cB_n)})$.   
But $ \gamma(W_i) = \frac 1 n V_i^*$ for all $1 \le i \le n$.
Thus $\overline{\gamma(\cB_n)}$ contains $n \ge 2$ co-isometries with 
mutually orthogonal initial spaces, 
and hence $\ltsr(\fA_n) = \ltsr(\cB_n)\ge \ltsr(\overline{\gamma(\cB_n)}) = \infty$.

\medbreak
(b) This proof is almost identical.   
Since the $\{ W_i \}_{i=1}^\infty$ are \emph{pure} co-isometries, the  \textsc{wot}-closed algebra 
$\fW_n$ generated by $\{ I, W_1, W_2, ..., W_n \}$ is a multiple of $\fL_n$, 
i.e. $\fW_n = (\fL_n^*)^{(k)}$ for some $1 \le k \le \infty$.  
Thus $\rtsr(\fW_n) = \ltsr ((\fL_n)^{(k)}) = \ltsr(\fL_n)$.   But the argument above used with the corresponding compression map
\[
 \begin{array}{rccc}
 \gamma: & \fW_n & \to & \bofh \\
		& X & \mapsto & P_\hilb X|_{\hilb}
 \end{array}	
\]
shows that $\rtsr(\fW_n) = \infty$, since $ \overline{\gamma(\fW_n)}$ contains at least $n \ge 2$ isometries with mutually orthogonal ranges.

Hence $\ltsr(\fL_n) = \infty$.   But by the Structure Theorem for free semigroup algebras mentioned above~\cite{DKP2001}, either there is
a homomorphism of $\fS_n$ onto $\fL_n$ or $\fS_n$ is a von Neumann algebra containing
two isometries with orthogonal ranges.  Either way, $\rtsr(\fS_n) = \ltsr(\fS_n) = \infty$.
\end{pf}



\bigskip

\bibliographystyle{plain}

\vskip 1 cm

\end{document}